# Majorization and the degree sequence of trees

[1]Leo Egghe and [2,3]Ronald Rousseau

[1]University of Hasselt, Belgium
leo.egghe@uhasselt.be
ORCID: 0000-0001-8419-2932

[2] University of Antwerp, Faculty of Social Sciences, B-2020 Antwerpen, Belgium, ronald.rousseau@uantwerpen.be

[3] Centre for R&D Monitoring (ECOOM) and Dept. MSI, KU Leuven, Leuven, Belgium, ronald.rousseau@kuleuven.be
ORCID: 0000-0002-3252-2538

**Abstract**

Purpose: We investigate the relation between degree sequences of trees and the majorization order.

Design/methodology/approach: We apply majorization using Muirhead's theorem.

Findings: We prove a theorem that provides a necessary and sufficient condition for delta sequences of trees to be comparable in the majorization order.

Research limitations: We only study trees, not general networks

Practical implications: Although our investigation is largely theoretical because trees are ubiquitous our study contributes to a better knowledge of trees as an important data structure.

Originality/value: This article is among the few combining Lorenz curves and majorization on the one hand, and degree sequences of networks on the other.

## 1. Introduction

In this introduction, we recall the notions used further in our article. These notions and their notation are well-known in network or graph theory, see e.g., [3], [9] or are taken from previous articles [1].

Let G = (V,E) be an undirected network, where $V = (v_k)_{k=1,...,N}$ denotes the set of nodes or vertices and E denotes the set of links or edges. We assume that #V = N > 1.

A path of length n is a sequence of vertices ($v_0$, …$v_k$, $v_{k+1}$, …, $v_n$) such that {$v_0$, …, $v_{n-1}$} and {$v_1$, …, $v_n$} are sets (being sets each consist of different elements) and for k= 0,…, n-1, $v_k$ is adjacent to $v_{k+1}$. A cycle is a path for which the starting point $v_0$ coincides with the endpoint $v_n$. A graph is connected if there exists (at least one) path between any two vertices. If #V = N, then the degree of node i, i = 1, …, N, i.e., the number of edges connected to node i, is denoted as $\delta_i$. In this article we always assume that G is connected, hence all degrees are strictly larger than zero. As there is no natural order among the nodes in a network we assume that these values are ranked in decreasing order.

Notation

The sequence of degrees of the nodes in a network G with N nodes is denoted as

$$\Delta_G = (\delta_1(G), \delta_2(G), ..., \delta_N(G)). \qquad (1)$$

We will informally refer to such a sequence as a delta sequence, consisting of delta values. Indices in the delta notation refer to a rank. Clearly, $\sum_{i=1}^{N} \delta_i = 2 (\#E)$, a notion which is known as the total degree of the network. It is easy to see that 2(N-1) ≤ $\sum_{i=1}^{N} \delta_i$ ≤ N(N-1). The lower bound is obtained e.g., for a tree (hence also for a chain) consisting of N nodes, while the upper bound is obtained for a complete graph where each node is connected to all other nodes.

Before moving on to examples and theory we recall the following definitions.

1.1 Definition: Trees, and branches

A free or unrooted tree is a connected graph with no cycles. Equivalently it is a connected graph such that removing any edge makes it disconnected. Another equivalent definition states that if v and v' are different vertices, then there exists exactly one path from v to v' [3].

Often there is one designated node, called the root. In that case, one says that the tree is rooted.

If *m* is any node in T (but not a terminal node, i.e. a node with degree one), then a branch rooted at *m*, consists of one link at *m*, and all nodes and links connected to *m* in T via that link. This is illustrated in Fig.1.

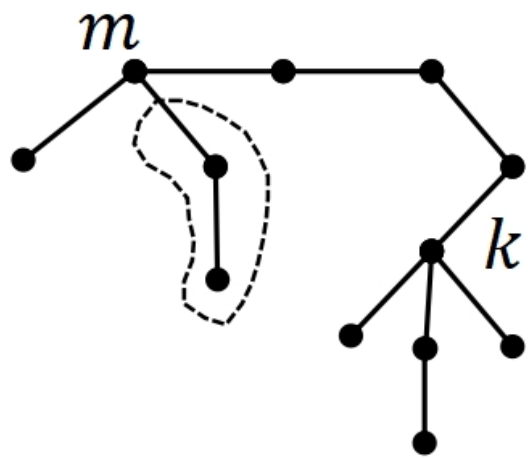

Fig.1. A tree and a branch rooted in node *m*

As emphasized by Knuth [3,p.305], trees are the most important nonlinear data structures.

1.2 Definition. Isomorphic graphs

Two graphs G and G' are isomorphic if there exists a bijection f between the vertices of G and G' such that there is an edge between vertices u and v in G if and only if there is an edge between the vertices f(u) and f(v) in G'.

When talking about a network or a tree we always mean the equivalence class of isomorphic networks or trees. Hence, we do not distinguish between isomorphic networks. By definition, two isomorphic networks have the same delta sequence, but the opposite is not true [1].

1.3 Definition: Spanning tree of a connected graph

A spanning tree of an N-node connected graph is a set of N-1 edges that connects all nodes of the network and contains no cycles. A graph may have different (non-isomorphic) spanning trees.

1.4 Definition: The Lorenz curve [4]

Let $X = (x_1, x_2, ..., x_N)$ be an N-sequence with $x_j \in \mathbb{R}^+, j = 1, ..., N$. If X is an N-sequence, ranked in decreasing order (always used in the sense that ranking is not necessarily strict), then the Lorenz curve of X is the curve in the plane obtained by the line segments connecting the origin (0,0) to the points $\left(\frac{k}{N}, \frac{\sum_{j=1}^{k} x_j}{\sum_{j=1}^{k} x'_j}\right)$, k= 1,…,N. For k = N, the endpoint (1,1) is reached.

1.5 Definition. The majorization property [5]

If X and X' are N-sequences, ranked in decreasing order, then X is majorized by X' (equivalently X' majorizes X), denoted as $X \leqslant_L X'$, if

$$\sum_{j=1}^{k} x_j \leq \sum_{j=1}^{k} x'_j \text{ for } k = 1, ..., N-1 \text{ and } \sum_{j=1}^{N} x_j = \sum_{j=1}^{N} x'_j \quad (2)$$

The index L in $X \preccurlyeq_L X'$ refers to the fact that this order relation corresponds to the order relation between the corresponding Lorenz curves. One may observe that X is majorized by X' ($X \preccurlyeq_L X'$) if and only if the Lorenz curve of X' is situated above (or coincides with) the Lorenz curve of X.

It is well-known, see e.g., [5, p.14] that $X \preccurlyeq_L Y$ is equivalent to each of the following statements:

(A) $\sum_i \varphi(x_i) \leq \sum_i \varphi(y_i)$ for all continuous, convex functions $\varphi: \mathbb{R} \to \mathbb{R}$.

(B) Y can be obtained from X by a finite number of elementary transfers [6].

Here an elementary transfer is a transformation from $(x_1, \cdots, x_N)$, where $(x_1, \cdots, x_N)$ is ranked in decreasing order, into $(x_1, \ldots, x_i + h, \ldots, x_j - h, \ldots, x_N)$ where $0 < h \leq x_j$.

1.6 Basic transfers

In the case that the elements in $(x_1, \cdots, x_N)$ are natural numbers, also h can be taken as a natural number, and it can even be taken to be equal to 1. In this case, we will say that this transfer is a basic transfer. It is shown in the appendix how to perform such basic transfers.

We write $X \prec_L Y$ for the strict Lorenz majorization, i.e., $X \preccurlyeq_L Y$ with X ≠ Y.

1.7 Definition: Non-normalized Lorenz curves

Let $X = (x_1, x_2, \ldots, x_N)$ be a decreasing N-sequence of non-negative real numbers, then the corresponding non-normalized Lorenz curve is the polygonal line connecting the origin (0,0) with the points $\left(j, \sum_{k=1}^{j} x_j\right)$, j = 1, …, N. This curve ends at the point with coordinates $\left(N, \sum_{k=1}^{N} x_j\right)$.

## 1.8 Definition: The non-normalized (or generalized) majorization order for N-sequences

If X and Y are decreasing N-sequences of non-negative real numbers, then X is majorized by Y, denoted as X ≼ Y if

$$\forall j, j = 1, \ldots, N: \sum_{k=1}^{j} x_k \leq \sum_{k=1}^{j} y_k \qquad (3)$$

The relation ≼ is only a partial order as non-normalized Lorenz curves (just like standard Lorenz curves) may intersect. If $x_j \leq y_j$, for $j = 1, \ldots, N$; then obviously X ≼ Y, but the opposite relation does not hold.

As for the Lorenz majorization, we write X ≺ Y, for X ≼ Y with X ≠ Y.

## 2. Basic transfers and delta sequences of trees

First, we explain the relation between a basic transfer and the delta sequence of a tree. Given, a tree T with delta sequence $\Delta_T = (\delta_1(T), \delta_2(T), \ldots, \delta_N(T))$, we know that always $\delta_N(T) = 1$. If now we perform a basic transfer, replacing the sequence $(\delta_1(T), \delta_2(T), \ldots, \delta_N(T))$ by $(\delta_1(T), \ldots, \delta_i(T) + 1, \ldots, \delta_j(T) - 1, \ldots, \delta_N(T))$, where $\delta_j(T) > 1$, we refer to this transfer as a basic tree transfer. Then we see that the degree of the node at rank j decreased by 1. This happens if we remove a branch (without the root node) from the node at rank j. As the degree of the node at rank i has increased by 1 and all other degrees have stayed invariant this can be realized by attaching a branch rooted at the node at rank j to the node at rank i This is illustrated in Fig.2 and in more detail in the appendix.

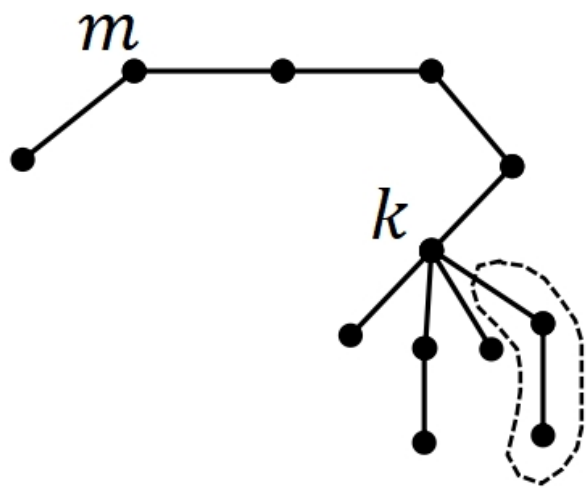

Fig.2. A branch at node *m* (see Fig.1) is replaced by the same branch placed at node *k*

3. Majorization and generalized majorization between delta sequences of networks

It is well-known that the relation $\preccurlyeq$ is not a total order between delta sequences of networks (of course with an equal number of nodes) (Egghe, 2024). We next show that this is not even true for trees. Consider the following five non-isomorphic trees with N = 8 nodes.

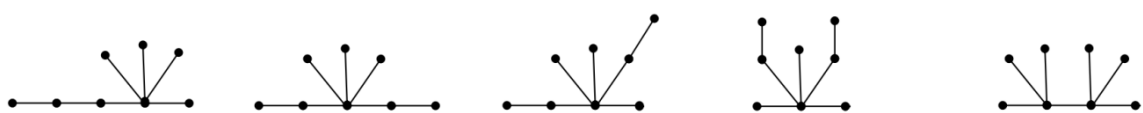

Fig.3. Five non-isomorphic trees with 8 nodes

The first four trees have the same delta sequence, namely $\Delta = (5,2,2,1,1,1,1,1)$ while the delta sequence of the last one is $\Delta' = (4,4,1,1,1,1,1,1)$. Clearly, $\Delta$ and $\Delta'$ are different and not comparable: $\Delta \not\preccurlyeq \Delta'$ and $\Delta' \not\preccurlyeq \Delta$.

Proposition 1

In the set of delta sequences of N-node trees, the majorization order coincides with the generalized majorization ($\leqslant = \leqslant_L$) and ($\prec = \prec_L$).

Proof. This follows immediately from the fact that the total degree of every tree with N nodes is 2(N-1).

Now we come to the main theorem of this article that provides a necessary and sufficient condition for delta sequences of trees to be comparable.

**Theorem**. Given two sequences $\Delta$ and $\Delta'$ of length N, where $\Delta$ is the delta sequence of a tree (hence does not contain a zero) and also $\Delta'$ does not contain a zero then

$$\Delta \prec \Delta'$$

$$\Leftrightarrow$$

For every tree T with delta sequence $\Delta$, there exists a tree T' with delta sequence $\Delta'$ which is created from the tree T by moving a finite number of branches to nodes with a higher or equal degree.

Proof. Assume that the tree T' is created from the tree T by moving a finite number of branches (each without their root) to a node with a higher or equal degree. If we replace in T with $\Delta_T = (\delta_1, \delta_2, \ldots, \delta_N)$ a branch of the node at rank j ($\delta_j > 1$) to the node at rank i where $\delta_i \geq \delta_j$ then only two values in $\Delta_T$ change (but we have no information about the new ranking): $\delta_i$ becomes $\delta_i + 1$ and $\delta_j$ becomes $\delta_j - 1$. The new delta sequence has values

$$(\delta_1, \ldots, \delta_{i-1}, \ \delta_i + 1, \ \delta_{i+1}, \ldots, \delta_{j-1}, \ \delta_j - 1, \delta_{j+1}, \ldots, \delta_N), \qquad (4)$$

perhaps in a different order. Anyway, the new delta sequence is strictly larger (in the $\prec\, =\, \prec_L$ ordering) than $\Delta_T$. Performing this operation a finite number of times proves that $\Delta \prec \Delta'$.

Conversely, we consider a tree T with delta sequence $\Delta$ (we know that such a tree exists). Hence the given sequence $\Delta$ is $\Delta_T$. We know now that $\Delta \prec \Delta'$. By Muirhead's theorem, we can apply a finite number of basic transfers on the tree T (moving from $\Delta$ to $\Delta'$). The resulting tree is the tree T' whose existence we have to prove. □

Figs. 1 and 2 illustrate this theorem with $\Delta = (4,3,2,2,2,2,1,1,1,1,1)$ and $\Delta' = (5,2,2,2,2,2,1,1,1,1,1)$.

The same reasoning as used in the Theorem can be used to prove the following well-known result.

Proposition 2 [2]. Given a sequence S of length N, consisting of strictly positive numbers, and with total degree 2(N-1), then we can construct a tree T such that S is the degree sequence of T, i.e. S = $\Delta_T$ .

Proof. If C denotes the degree sequence of the N-node chain, then $C \prec_L S$, using that they both have a total degree of 2(N-1). Now, by Muirhead's theorem, we can apply a finite number of basic transfers to C and reach S. The resulting tree is the tree T whose existence we have to show.

Remarks

Remark 1. Proposition 2 is Corollary 1, p. 499 in [2].

Remark 2. The obtained tree T' does not have to be unique as illustrated in the appendix.

Remark 3. The theorem states "for all T with delta sequence Δ, there exists a tree T' with delta sequence Δ' ". The theorem is false when this expression is replaced by "for all T and T' with delta sequences $\Delta_T = \Delta$ and $\Delta_{T'} = \Delta'$". We provide an example for N=8. Let T be the tree shown in Fig.4 (a) and T' the tree shown in Fig. 4 (b). Then $\Delta = \Delta_T$ = (4,2,2,2,1,1,1,1), $\Delta' = \Delta_{T'}$ = (5,2,2,1,1,1,1,1) and $\Delta \leqslant \Delta'$. Yet, it is impossible to transform T into T' via basic transformations (recall the condition $\delta_j \geq \delta_i$).

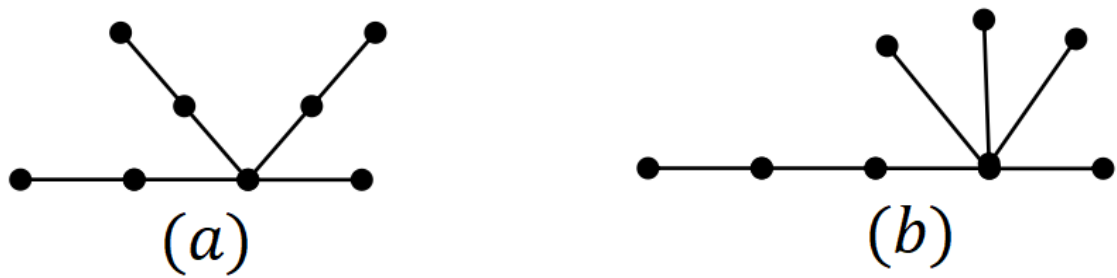

Fig.4 Trees T and T'

Using this example we can construct an illustration of the theorem though. Consider the tree T" in Fig. 5.

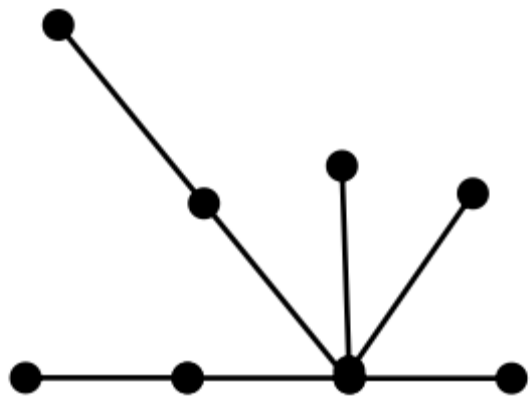

Fig. 5 Tree T"

Then $\Delta" = \Delta_{T"}$= (5,2,2,1,1,1,1,1) and T" can be obtained from T by basic transfers.

Remark 4. The relation $\preccurlyeq = \preccurlyeq_L$ is a total order for trees if and only the number of nodes N ≤ 7.

Proof. All non-isomorphic trees for N < 23 can be found at https://users.cecs.anu.edu.au/~bdm/data/trees.html. Then one can check that for N ≤ 7 we have a total order. For N=8 we already gave an example that the order is not total. Based on this example it is easy to construct examples for all N > 8, see Fig. 6. This figure gives a tree T (left) with $\Delta_T = (5, \underbrace{2, \ldots, 2}_{N-6\ times}, 1,1,1,1,1)$ and on the right, a tree T' with $\Delta_{T'} = (4,4, \underbrace{2, \ldots, 2}_{N-8\ times}, 1,1,1,1,1,1)$. Then $\Delta_T \not\preccurlyeq \Delta_{T'}$ and $\Delta_{T'} \not\preccurlyeq \Delta_T$.

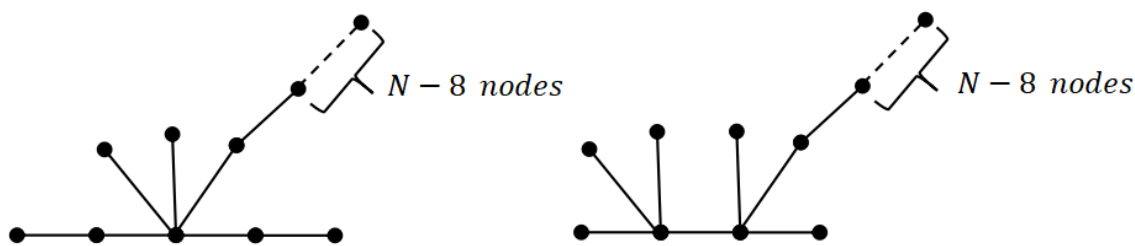

Fig. 6. Non-isomorphic, non-comparable trees for N > 8

Proposition.

If $\Delta_C$ is the delta sequence of the N-node chain and $\Delta_M$ is the delta sequence of a connected N-node network M (not being a chain), then

$$\Delta_C \prec \Delta_M$$

Proof. Let $\Delta_T$ be the delta sequence of any tree T. Then we know that we can obtain this tree T be a finite number of basic transfers from a chain and hence by the theorem $\Delta_C \preccurlyeq \Delta_T$, with equality only if T is a chain. Consider now any N-node network M, then this network has a spanning tree $T_M$ with $T_M \preccurlyeq \Delta_M$. By the transitivity of $\preccurlyeq$ and the fact that M is not a chain, we obtain that $\Delta_C \prec \Delta_M$.

4. Applications

In data file systems, directories and files are often represented as a tree structure. Moving a subdirectory (and all its contents) from one directory to another is a common operation. If the target directory has a higher degree than the original one, this means, in the terminology of this article, that the new situation majorizes the old one.

Suggestions to restructuring the fossil record, (phylogenetic trees), leading to changes akin to moving branches in a tree happen regularly, see e.g., [8].

5. Conclusion

This article is among the few combining Lorenz curves and majorization on the one hand, and degree sequences of networks on the other. It supports Rousseau's statement [7] that Lorenz curves (and hence the majorization order) are universal tools for studying networks.

In particular, we proved a theorem that provides a necessary and sufficient condition for delta sequences of trees to be comparable in the majorization order. Our methodology leads to an almost trivial proof of Hakimi's corollary (a corollary of a more general result about linear networks) on the realizability of a set of strictly positive natural numbers as degrees of the vertices of a tree.

Because trees are ubiquitous our study contributes to a better knowledge of trees as an important data structure.

**Acknowledgments.** The author thanks Li Li (Beijing, China) for drawing excellent illustrations and Raf Guns (Antwerp, Belgium) for useful discussions.

The authors declare no conflict of interest. No funding was received for this work.

**References**

[1] L. Egghe. Networks and their degree distribution, leading to a new concept of small worlds. *Journal of Informetrics* (2014; to appear).

[2] S.L. Hakimi. On realizability of a set of integers as degrees of the vertices of a linear graph I. *Journal of the Society for Industrial and Applied Mathematics,* **19**(3) (1962) 496-507.

[3] D.E. Knuth. *The Art of Computer Programming (second ed.) Vol.1 Fundamental Algorithms.* Reading: Addison-Wesley, 1973.

[4] M.O. Lorenz. Methods of measuring the concentration of wealth. *Publications of the American Statistical Association*, **9** (1905), 209-219.

[5] A.W. Marshall, I. Olkin and B.C. Arnold. *Inequalities: Theory of Majorization and its Applications*. New York: Springer, 2011.

[6] R.F. Muirhead. Some methods applicable to identities and inequalities of symmetric algebraic functions of n letters. *Proceedings of the Edinburgh Mathematical Society*, **21** (1903) 144-157.

[7] R. Rousseau. Lorenz curves determine partial orders for comparing network structures. *DESIDOC Journal of Library & Information Technology*, **31**(5), (2011) 340-347.

[8] E. Tschopp, O. Mateus, and R.B.J. Benson. A specimen-level phylogenetic analysis and taxonomic revision of Diplodocidae (Dinosauria, Sauropoda). *PeerJ* **3**, (2015) e857. https://doi.org/10.7717/peerj.857

[9] S. Wasserman and K. Faust. *Social Network Analysis*. Cambridge (UK): Cambridge University Press, 1994.

**Appendix**

We provide an algorithm in pseudo-code to transform $Y = (y_1, \cdots, y_N)$ to $X = (x_1, \cdots, x_N)$ with $Y \leqslant_L X$ by performing basic transfers. We assume that X and Y are ranked in decreasing order and that they are not equal (otherwise nothing must be done). As X and Y are trees we know that $x_1$ and $y_1$ are both strictly larger than 1 (for N > 2) and $x_N$ and $y_N$ are equal to 1.

For i = 1 to N-1  (i represents an index)

    While $y_i$ < $x_i$

    Find j such that $y_j$ > $x_j$ (otherwise the transfer cannot lead to the required result)

    Apply a basic transfer (transfer by 1) from node $y_j$ to node $y_i$

    Reorder Y

If Y = X the algorithm ends.

Example. For N = 8: Y = (3,3,3,1,1,1,1,1) ⩽ X = (5,3,1,1,1,1,1,1)

We take i = 1 and observe that 3 < 5. Next we see that $y_3$ = 3 > $x_3$ = 1, hence j = 3. We apply a basic transfer leading to Y = (4,3,2,1,1,1,1,1)

Still with i = 1, (as 4 < 5), we have j = 3, with $y_3 = 2 > x_3 = 1$.

Again we apply a basic transfer leading to Y = (5,3,1,1,1,1,1,1) = X

For the corresponding trees, we have (for example):

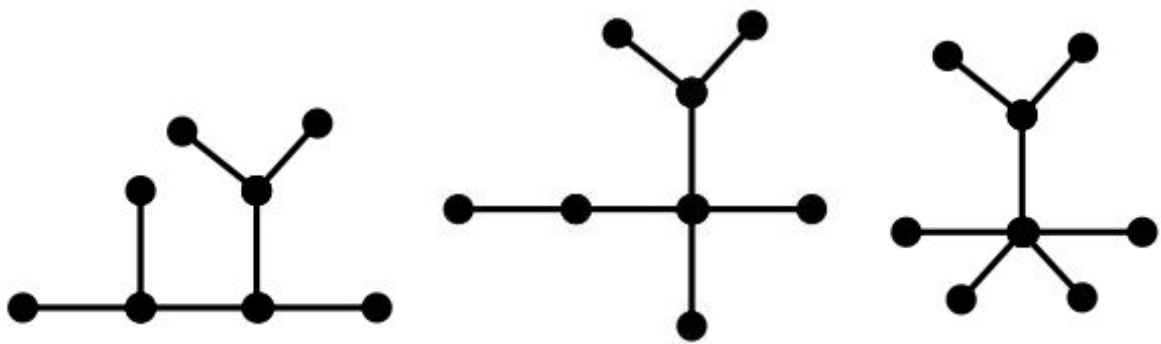

Fig. 7. Moving branches to go from Y(left) to X(right)

Next, we provide an example where we start from the 8-node chain leading to Δ = (5,2,2,1,1,1,1,1), the delta sequence of the first four trees of Fig.3. Recall that the delta sequence of the 8-node chain is (2,2,2,2,2,2,1,1)

Following the algorithm, we see that, for i =1, 2 < 5, hence j = 4 (as 2 > 1).

This leads to (3,2,2,1,2,2,1,1) and rearranging gives: (3,2,2,2,2,1,1,1).

Now, i in the algorithm is still equal to 1 (3 < 5) and j = 4 (2 > 1) leading to (4,2,2,1,2,1,1,1). Rearranging gives: (4,2,2,2,1,1,1,1).

Still i =1 (4 < 5), and j=4 (2 > 1) which leads to (5,2,2,1,1,1,1,1) = Δ.

The corresponding trees are not unique as they depend on the indexing of the nodes. The next figure shows how to get from the chain to a tree with delta sequence Δ. The nodes are indicated by their index number (possibly changing in each step). Observe that for two nodes with equal degrees, we are free to index them as we wish.

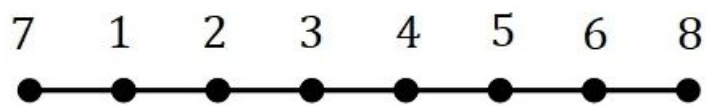

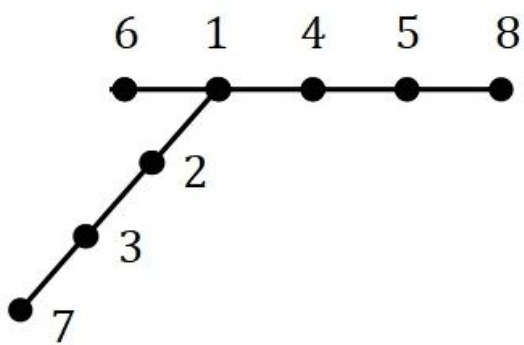

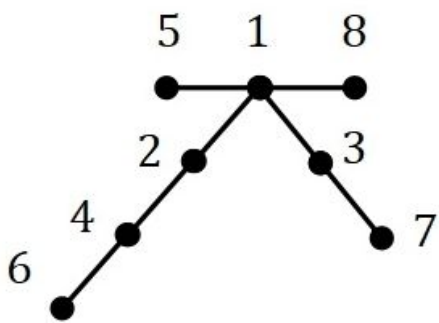

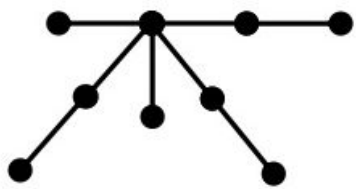

Fig. 8. Algorithm leading to the 4<sup>th</sup> tree in Fig.3

Next, we apply a different way of indexing the nodes. This leads to the third tree of Fig. 3.

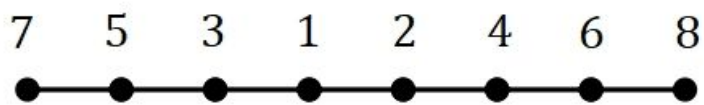

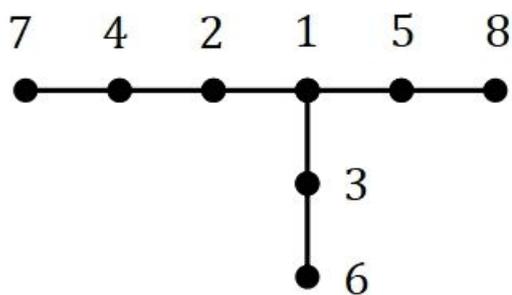

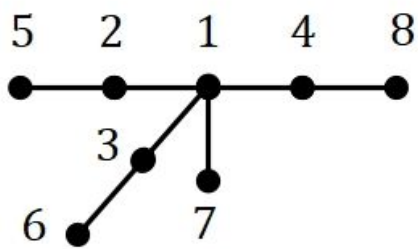

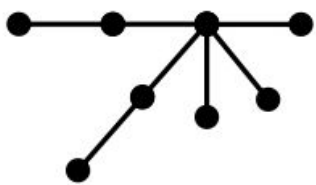

Fig. 9. Algorithm leading to the 3<sup>rd</sup> tree in Fig.3.